# A NOVEL METHOD FOR PREDICTION AND APPROXIMATION OF FUNCTIONS (SELF APPROXIMATION METHOD)

M.Abolqhasemi, F.DIDEHVAR, E.SAFAVIYEH, N.HASHEMI

1. Introduction

The natural laws, which can be recognized by mind based on induction, govern the nature in an almost periodic manner. As regular laws in nature, the change of season and the climate diagram of the temperature have such an almost periodic behavior.

It is good to remind that there are many similar functions in different areas, such as

Signal processing, vision, economy, etc.

In general, such regular functions can be written as a summation of some basic periodic functions as basis. The two main examples which support this concept are Fourier series and wavelets. As we know, the using these methods for approximations have limited applications in prediction.

One of the alternative ideas which is proposed in this paper is to consider function f itself in order to find basic functions to approximate f, instead of sine and cosine functions in Fourier series or other basic functions in wavelet methods. Based on this idea and finding basic functions which are derived by f itself, as we will show it later, this method introduces a better performance for predicting the future of some empirical data.

Generally, previous methods like Fourier series and wavelets are not considered to solve forecasting problems, but by means of such new method which is introduced in this paper, we can obtain an extension of these methods to propose a new approach in the area of prediction. Therefore, our

___________________________________________________________________________________

The names are in alphabetic order. This article is written based on a project in problem solving group (khawrazmi group).

method is motivated by an analogy to the main idea which is applied by Fourier and wavelets series procedure.

In order to implement this method we first define the notion of "semi-period" of a function which is the generalized concept of the period of a function. Despite the periodic functions, the "semi-periodic" functions are defined for the wide variety of applications, although they don't seem periodic.

**Our New method**

The methods of Fourier series and wavelets have been applied for a long time in order to approximate functions. But in the scope of prediction they have limited usage.

The central idea behind theses methods is to build some basic functions like sine and cosine in Fourier series and modeling the behavior of a function by a weighted sum of these basic functions.

In this article, the same idea has been applied but instead of sines and cosines and wavelets we use the original function itself in order to produce new basic functions.

This new method intuitively supports the idea that the function itself has more information for constructing and predicting a function rather than applying sines, cosines, or wavelets.

Considering the function itself for producing the basic functions, we could provide a method to predict the path of a function following a regular rule. Like the previous methods, we apply this method for approximating a function by using the basic functions, derived from the original function itself.

**The Procedure of the New Algorithm**

**Step 1:** How to find a **semi period**.

We know that a function is predictable, if we can find a regular rule in the behavior of the function. This rule could be defined in terms of a period. Although most of the functions that we know are not periodic, there is a large class of functions which could be considered (intuitively) **almost periodic**. For example, the function, shown in Fig 1, is not periodic at first step, but we can consider them intuitively as semi-periodic functions. The concept of a period can be extended to a **semi period** which is used throughout this paper. In the following, we will give a definition for the concept of a **semi period**.

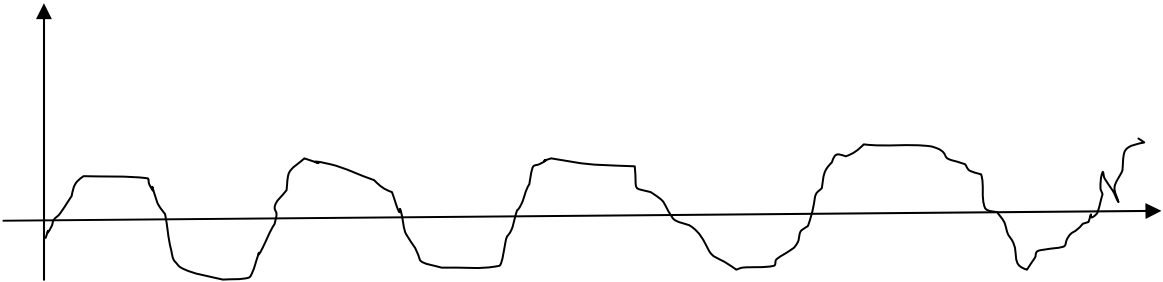

Let f(x), 0 <= x <= t, represent the real data measured in the interval [0,t]. Suppose T is a period, and S(n) is defined as a partitioning function which divides the function f into a number of equally spaced segments {F1,F2,…,Fn}.

Obtained by partitioning the function f with the length of T in the interval [0,t], all segments are equally spaced with the length T, except the last segment which has a length less than t.

All segments have lengths T except the last segment which has a length less than (t?).

Suppose that S is a partition of the interval [0, t] to equal segments, each segment has length T except the last segment which is equal or less than T. Now we define the partition Fi as follows

$$F_i(x) = f(x + (i-1)*T+1), \text{ for } 0<x<T, i = 1, 2, L \text{ such that } k = 1,…, n$$

where T is the length of the initial partition. (If we assume that the function f is the measured data points at discrete times (time?) k, then we introduce the following formula). …

.

Suppose that X1...Xm is the instances in our discrete space. Now we give a rearrangement of naming these instances as following

$$LSG(T) = \sum_{i=1}^{n} \frac{\sum_{j} \frac{|F_{i+1}(x_j) - F_i(x_j)|}{L_i} \times \frac{L_i}{T}}{n + |F_{n+1}|/T},$$

$$L_i = length(F_i)$$

**LSG FORMULA (1)**

Now we define a semi period as the minimum of the series given in (1).

Suppose that the time interval T* is defined as the semi period of the function f. Clearly, for a periodic function, the value of T* will be exactly equal the period of the periodic function.

These formulas can be applied to the data of time series in naturalistic events to find a semi period T* which is almost familiar to everyone, like the number of days in a week or, in a larger scale, the number of years, etc.

We can apply this semi period concept to forecast functions containing some kind of regularities. The procedure of this forecasting follows the algorithm in the next steps.

**The next steps of the algorithm**

The periodic functions are oscillating around a horizontal line. We expect that the semi period functions also oscillate approximately around such a horizontal line. But the trend of many functions which we know, does not satisfy this constraint. Therefore, in the next steps we follow not only the concept of semi period but also we will consider the trend of a function.

In a general view, we first estimate the trend of function as Tf and then define the error (or noise) function G = f – Tf. Now the function G is a better choice for finding its semi period rather than using the original function f, because G oscillates around a line parallel to the X axis. It is important to note that our method is very sensible to the trend of the function and we, therefore, need to find a good approximation for this trend. In the following we will see the procedure of this method.

**Step2& Step 3: estimate series with basic elements and series construction with basic elements**

After finding quasi-period out from LSG formula, basic elements could be constructed to apply for series of estimation series. Roughly speaking, the procedure that makes these basic elements is to divide
The interval into subintervals in the length of the quasi period, and considering the sub functions
Which are made, transferring all to the first interval, and then we get the average.
After finding quasi-period and retrieving basic elements, we intend to construct series in the third step.
In this step, a simple rule is used. This rule is governed by primitive idea that repetition phenomena can be known by only one repetition pattern of phenomena.
Taking into account this primitive rule implies that series construction is possible with repetition of basic element. Therefore, new series is a series of basic elements sequential that repeated during at its interval.
Mathematically, the result would be the using of basic element in order to produce an approximation of the original function. This is a pattern which by drawing it in all intervals we find out the first approximation of the function.

As you will see, because of using quasi-period value for constructing basic element, our main goal to use series itself could be satisfied. An overview over of this method shows basic elements which are constructed according to above procedure can meet a wide range of criterion, as mentioned in previous description.
Under such assumption, estimating functions with this basic element could be applied for more cases other than merely representing series in some order that Fourier and wavelet do. This could lead us to construct series out of given data, not only approximating functions, this means that we can use it as a forecasting technique. The power and weakness of this technique will be presented in the examples at the end of the paper.

The mechanism of how basic elements arrange, makes new method for constructing series. For example applying weighting element according to position provides a new approach in constructing series. Also, in this section for constructing series a wide range of intelligent algorithms could be applied in order to determine weights of best fit. A number of these methods are applied in examples in that would be revealed in the following sections.

**Step4: The subject matter of a trend**

As usual, any natural series that have been designed in all prediction algorithms to forecast contain some trends. In particular, the time series can be considered as a composition of deterministic and nondeterministic components modulated through a trend. So, the nice distinction between such a trend and the other components can be discussed.
Suppose the series constructed by considering daily revenue of a normal shop for a period of ten year. It can be seen that the shape of figure that such time series has is altering increasingly. There are some short term rises and falls in it but in general it is growing. This growth might be affected by inflation in prices. In this case, the rate of inflation is so wide that makes it difficult to estimate its underlying behavior. Perhaps if there were some more data, e.g. 100 years, there had been a chance to predict the behavior of inflation, but as far as we are concerned we have no solution but model it by using a trend concept.

Although there are a lot of traditional methods for de-trending, such as differentiating or taking a weighted moving average, we have simply demodulated the trend with polynomial functions, because it has more consistency with other parts of our method. In general our method is to do a regression with polynomials of degrees between 0 and max degree and then subtracting the values of the regression from the values of the series. The max degree is a definite value defined by the user, We set it here, on 10. Among the regressions with different polynomial degrees the most beneficial of it one must be selected. This selection is made based on validation set that has been explained in following sections in selection procedure.

We also consider an innovative method to take the fractional degrees of polynomials into account. Suppose that a polynomial of degree *m* is the most perfect one.(f-m) and a polynomial of degree m+1is our choice for the polynomials of this degree.( f-m +1).It is reasonable to try to test some functions between polynomials of degree *m* to *m+1,* which they fit more convenient for the trend. We use the following equations to consider such functions between polynomials of degree m and m+1 in general.

$$f_{m,N}^{i}(x) = \sqrt[N]{f_m(x)^i f_{m+1}(x)^{N-i}}$$

By the above equation The interval of polynomials of degrees *m* and *m+1* is divided into *N* functions and the $f_{m,N}^{i}(x)$ is the $i^{th}$ one. So we have N function here. In conclusion the most beneficial one among these functions on validation set can be considered as the best.

**Step5: Evaluation and selection**
In general our algorithm, at each iteration, extracts some information from the given data and reforms its model of the underlying process. This extracted information is in the form of a periodic series. As it was already mentioned there are some factors involved in estimating this periodic series. The examples of which are methods of approximating a period and the degree of polynomial used for de-trending. So, at each iteration, there are a lot of factors and above all a set of configurations and one should select the best one among them.

In order to select the best configuration, we divided the given data into three sets. One set, the largest one, is called train set. The train set is the set which is, as we assume given completely. All parameter estimations such as the period of series must be done according to it and no other data should be given by outsourcing.

The other set is called the test set, it is the one by which we test our algorithm in order to show its performance. This set must be considered as one that is not given. In machine learning, if a learning algorithm could have got some good results in this set, it had been a general and perfect algorithm.

As it is proposed we should select some factors among the configurations. One may say that it is possible to select the best configuration by dealing their performance on test set and select the best one, but it is an unfair strategy, because the test set must be regarded as an undisclosed one until the end of algorithm and revealing data from it to select the best configuration is wrong. So another set may be needed to pick out the most satisfactory parameters. We use a set called validation set.

Validation set is that which is not given, as assumed, through parameter estimation but is given only for model selection. So, we make a set from different configurations which are set through training phases of each iteration, and then select the best one through validation phase and repeat this procedure for next iterations. At the end of algorithm operation we can test the final model on test set to show its generality and performance.

The performance out puts can be measured through different error functions like MSE, NMSE, MAE, MAPE, and SMAPE. We use SMAPE or Symmetric Mean Absolute Percentage Error in our algorithm because it can be computed through equation #:

$$SMAPE(Y,F) = \frac{1}{n}\sum_{t=1}^{n}\frac{|Y_t - F_t|}{(Y_t + F_t)/2} *100$$

Where *F* is the actual values of the series and *Y* is the approximated one.

The size of training and validation sets can also be varied and this change may provide a change in the performance of prediction. So, we have used different ratios of training and validation size. This ratio should also be adopted; therefore, it must be included in the selection procedure. In order to make a fair selection the different configurations must be evaluated in a common dedicated part of the validation sets. Thus, a configuration is the best if it has the best performance.

**The complete explanation of the algorithm**

In order to get a good vision of underlying logic by reviewing precedent sections, we tried to explain the algorithm in its developmental cycle. In this section the final algorithm is proposed. Although the algorithm is simple in general, it has complicated details, so it should not have been proposed directly.

In general the algorithm is to detect a periodic behavior in data, model it and then do the same for the rest. In more details it does the following for each iteration:

1- De-trending
2- Find the period length (like-period)
3- Create a series in a length of time (basic-element)
4- Repeat the basic element in order to estimate of this iteration
5- Subtract the series of last step from the series of this iteration and construct the series of next iteration

The algorithm may end with some stopping criteria like a maximum iteration, or meet a worse situation for some iteration.

There are a lot of variations of algorithm to do so, various degrees for polynomials to detrend, different methods to find the period, some weighting methods to create the mean basic-element series, and some different validation and train sizes. To take these varieties into account the algorithm must be more complicated as detailed more in the following items:

- ❖ In case stopping criteria is not met:
- ➢ For various validation ratios:
- ▪ Get the train part of the series and code it as f{i}
- ▪ For different period finding methods:
- ● For some weighting strategies to calculate mean series:
- ♦ For different degrees of polynomials:
- ➢ De-trend the f{i}. i.e. f{i}=f{i}-REG (regression series)
- ➢ Find the like-period of f{i}
- ➢ Calculate the mean basic-element series by using weighting method
- ➢ G=Repeat the basic-element series (side by side) to make an estimate of f{i}
- ➢ finalPredict = finalPredict + REG + G
- ➢ calculate the validation error
- ♦ select the best one of above and do the same for fractional polynomials near the best
- ▪ select the best of above and do the same for other validation sizes
- ➢ select the best parameters of all and let f{i+1}=f{i}-G (corresponding to the best configuration)
- ❖ end

**Results and analysis**
**Good results**

| Series name | Train & validation period | Test period | Train NMSE | Train SMAPE | Test NMSE | Test SMAPE |
|---|---|---|---|---|---|---|
| HD[1] | 3HD | | 0.033823 | 1.6043 | 1.0436 | 8.0042 |
| HD | 4HD | | 0.043507 | 1.4857 | 1.2919 | 9.0334 |
| HD | 5HD | | 0.016168 | 0.45198 | 0.20144 | 3.351 |
| HD | 6HD | | 0.020163 | 0.48124 | 0.3529 | 5.745 |
| HD | 7HD | | 0.0060526 | 0.23152 | 0.30834 | 5.0753 |
| Gold | 10Gold | | 0.021122 | 0.21009 | 4.2312 | 1.5508 |
| Gold | 1Gold | | 0.0052454 | 0.26024 | 2.6908 | 7.7231 |

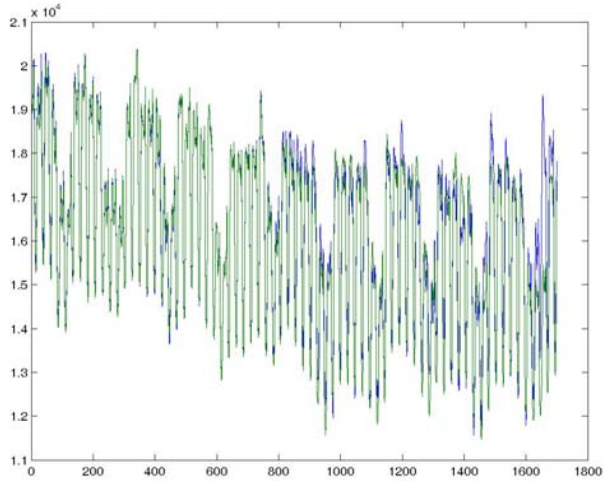

3HD(Hour, Demand 11-mar-08 13:00 to 21-may-08 8:00)

---

[1] Hourly Demand (Ontario Energy records)

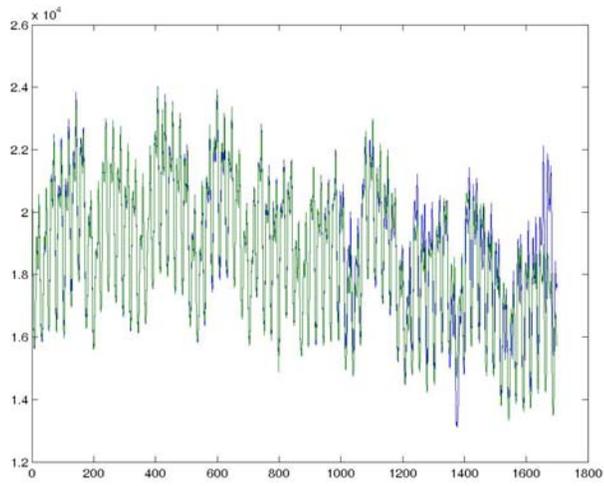

4HD (Hour, Demand   19-jan-07 21:00 to 31-mar-07 16:00)

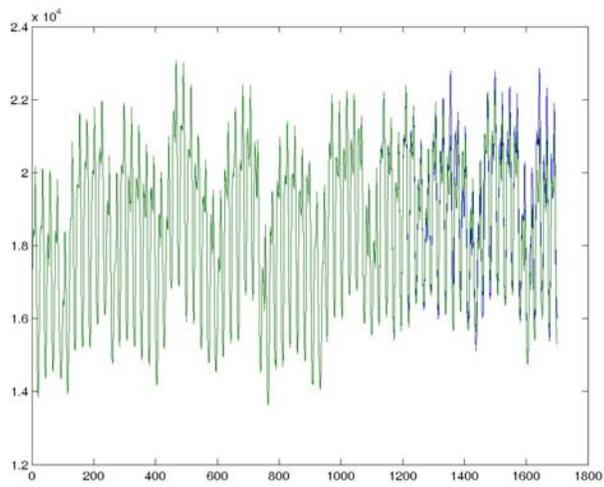

5HD (Hour, Demand 28-dec-05 9:00 to 11-jan-06 5:00)

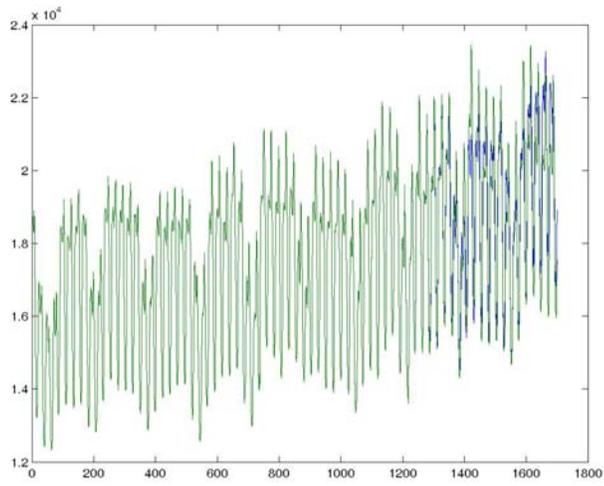

6HD (Hour, Demand 8-oct-04 13:00 to 7-nov-04 21:00)

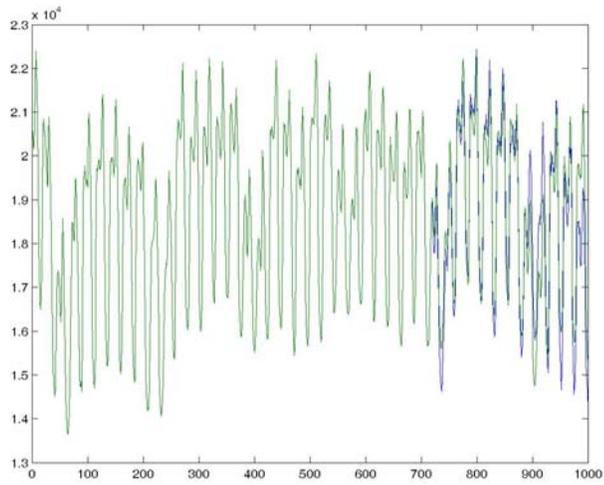

7HD (Hour, Demand 26-jan-06 13:00 to 9-mar-06 4:00)

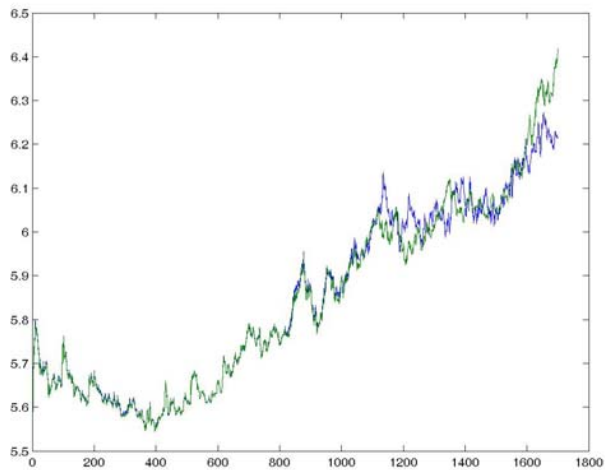

10Gold(Day, Price)

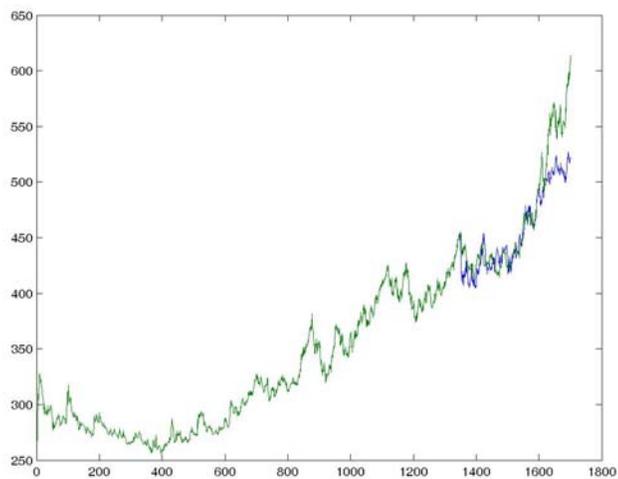

11Gold (Day, Price)

**Middle cases**

| Series name | Train& validation period | Test period | Train NMSE | Train SMAPE | Test NMSE | Test SMAPE |
| --- | --- | --- | --- | --- | --- | --- |
| HOEP | 9HOEP | | 0.1753 | 13.505 | 1.698 | 32.235 |
| Gold | 11Gold | | 0.024812 | 0.07856 | 34.445 | 2.3891 |

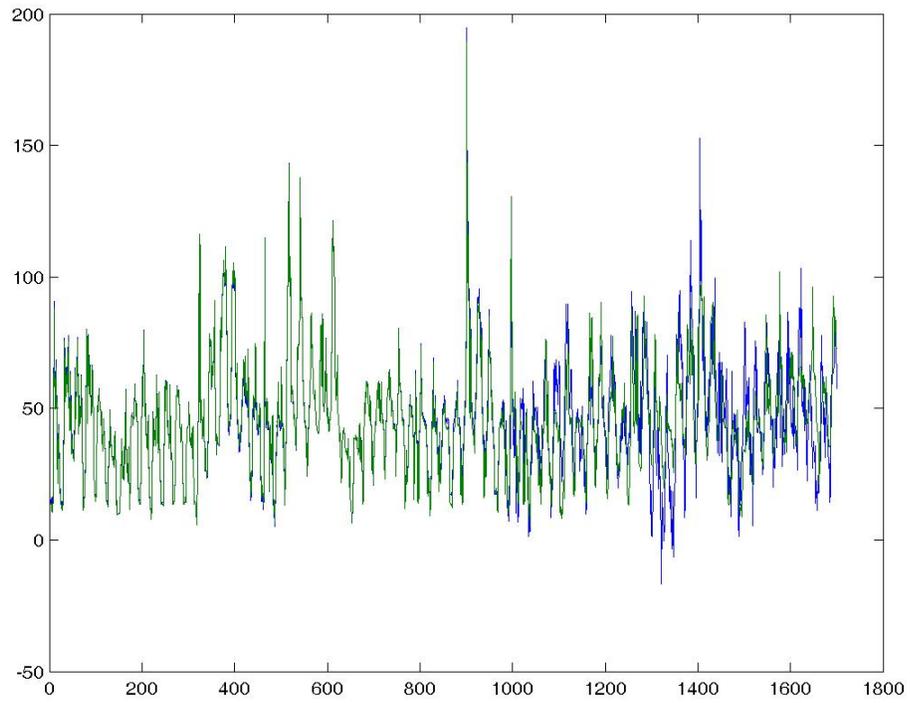

9HOEP(Hour, Price 24-may-04 1:00 to 2-aug-04 20:00)

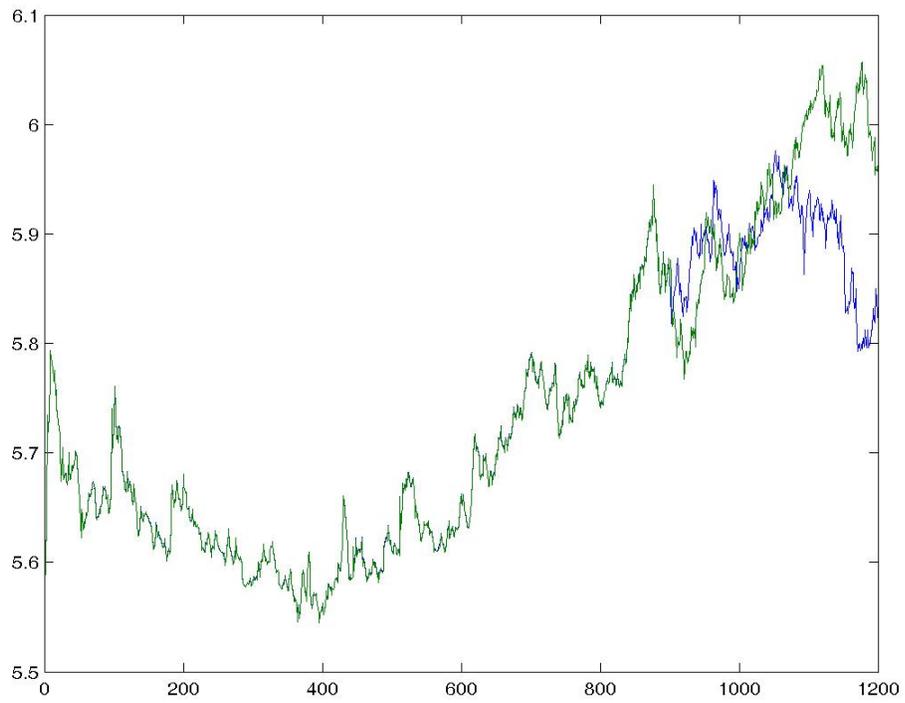

11Gold(Day, Price)

**Worse cases**

| Series name | Train& validation period | Test period | Train NMSE | Train SMAPE | Test NMSE | Test SMAPE |
|---|---|---|---|---|---|---|
| HOEP | 8HOEP | | 0.29795 | 10.961 | 1.3048 | 38.546 |

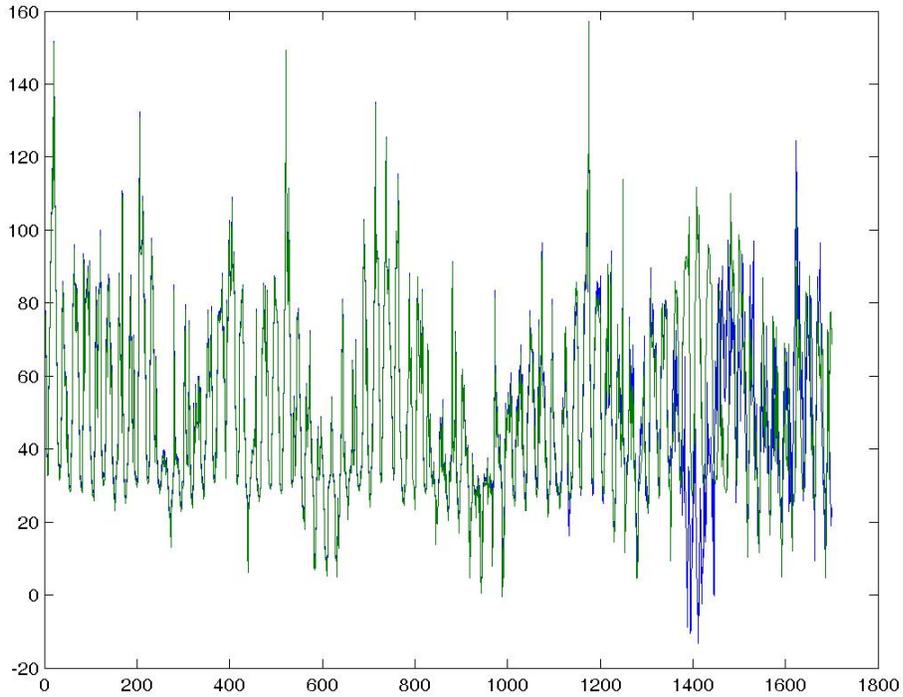

8HOEP ( Hour, Price 7-aug-07 21:00 to 17-oct-07 16:00)

**Conclusion**

Throughout this article the major conclusion is about comparing this method with some very famous methods like Fourier series and wavelets. As it was mentioned before, this method provides a way to apply function itself instead of some prior known functions such as sine and cosine functions in Fourier series, or wavelets. In brief, the power of interpolation and approximation of a function by this technique is as powerful as, we could use it efficiently in extrapolations and predictions in a wide range of examples, while we couldn't say the same for Fourier series and wavelets.

These are the central ideas in this article, the foundation this article is based on.

There is only one problem here, about interpolation. In fact this method, in most practical purposes
give a satisfactory approximation of a function in a compact interval, but mathematically there is no proof for that, on the contrary, it seems in the concluded process of "approximation method" the residue didn't tend to zero.

In modifying this method, we could present an algorithm which satisfies the above theoretical purpose for approximation. In fact, if we add a condition in the part of algorithm which is specified to find the trend function , which based on that we choose our trend function such that it gives us an approximation of the function better than %50 (in each iteration), which satisfies our desired goal in approximation. It is sufficient to note that all the time to do this approximation, by considering higher orders polynomials is possible. But nevertheless this modification could violate some practical purposes for predictions (as we have checked it in some examples).

We know this method somehow a deterministic method. Applying and different combination of this method and some stochastic methods could make powerful methods for predictions.